\documentclass[a4paper,11pt]{amsart}
\usepackage{amssymb,amsxtra,eucal}
\usepackage{fullpage}
\usepackage[all]{xy}\CompileMatrices

\theoremstyle{plain}
\newtheorem{thm}{Theorem}[section]
\newtheorem{lem}[thm]{Lemma}
\newtheorem{prop}[thm]{Proposition}
\newtheorem{cor}[thm]{Corollary}

\theoremstyle{definition}
\newtheorem{rem}[thm]{Remark}

\newtheorem{ex}[thm]{Example}

\numberwithin{equation}{section}

\def\eg{\emph{e.g.}}
\def\ie{\emph{i.e.}}
\def\ds{\displaystyle}
\def\:{\colon}
\def\.{\cdot}

\def\<{\left\langle}
\def\>{\right\rangle}
\def\({\left(}
\def\){\right)}
\def\ph#1{\phantom{#1}}
\def\epsilon{\varepsilon}
\def\phi{\varphi}
\def\subset{\subseteq}

\def\leq{\leqslant}
\def\geq{\geqslant}

\def\lra{\longrightarrow}
\def\Lra{\Longrightarrow}
\def\ra{\rightarrow}
\def\hat#1{\widehat{#1}}

\def\iso{\cong}

\DeclareMathOperator{\im}{im}

\def\F{\mathbb{F}}

\def\Q{\mathbb{Q}}

\def\Z{\mathbb{Z}}

\def\Times_#1{\mathop{\times}_{#1}}
\def\oTimes_#1{\mathop{\otimes}_{#1}}
\def\oPlus_#1{\mathop{\bigoplus}_{#1}}
\def\ideal{\triangleleft}

\DeclareMathOperator{\Cotor}{Cotor}

\DeclareMathOperator{\Tor}{Tor}

\def\id{\mathrm{id}}

\DeclareMathOperator{\sign}{sign}
\DeclareMathOperator{\incl}{incl}
\DeclareMathOperator{\cone}{cone}
\def\cp{\textsl{cp}}
\def\StA{\mathcal{A}}
\def\StB{\mathcal{B}}
\def\StE{\mathcal{E}}

\begin{document}
\title[Cooperations for the connective Adams summand]
{On the cooperation algebra of the connective Adams summand}
\author{Andrew Baker \and Birgit Richter}
\address{Mathematics Department, University of Glasgow,
Glasgow G12 8QW, Scotland.}
\email{a.baker@maths.gla.ac.uk}
\urladdr{http://www.maths.gla.ac.uk/$\sim$ajb}
\address{Fachbereich Mathematik der Universit\"at Hamburg,
20146 Hamburg, Germany.}
\email{richter@math.uni-hamburg.de}
\urladdr{http://www.math.uni-hamburg.de/home/richter/}
\thanks{
We would like to thank Iain Gordon, John Rognes and Sarah Whitehouse for
helpful comments. The first author was supported by the Max-Planck institute
of mathematics, Bonn, and  the second author  by the
\emph{Strategisk Universitetsprogram i Ren Matematikk}
(SUPREMA) of the Norwegian Research Council.
We thank the universities of Bern, Bonn, and Oslo for their hospitality.}
\keywords{connective $K$-theory}
\subjclass[2000]{Primary 55P43, 55N15; Secondary 55N20, 18G15}
\begin{abstract}
The aim of this paper is to gain explicit information about the multiplicative
structure of $\ell_*\ell$, where $\ell$ is the connective Adams summand.
Our approach differs from Kane's or Lellmann's because our main technical
tool is the $MU$-based K\"unneth spectral sequence. We prove that the algebra
structure on $\ell_*\ell$ is inherited from the multiplication on a Koszul
resolution of $\ell_*BP$.
\end{abstract}
\maketitle

\section*{Introduction}

Our goal in these notes is to shed light on the structure, in particular
on the multiplicative structure, of $\ell_*\ell$, where we work at an odd
prime~$p$ and~$\ell$ is the Adams summand of the $p$-localization of the
connective $K$-theory spectrum~$ku$. This was investigated by
Kane~\cite{Kane:AMSMem} and Lellmann~\cite{Lellmann:Ops&Coops} using
Brown-Gitler spectra. Our approach is different and exploits the fact that
$MU$ is a commutative $\mathbb{S}$-algebra in the sense of Elmendorf, Kriz,
Mandell and May~\cite{EKMM} and $\ell$ is a $MU$-ring spectrum (in fact it
is even an $MU$-algebra). As a calculational tool, we make use of a K\"unneth
spectral sequence~\eqref{eqn:KunnSS-l*l} converging to $\ell_*\ell$ where we
work with a concrete Koszul resolution. Our approach bears some similarities
to old work of Landweber~\cite{Landweber:MU*ku}, who worked without the
benefit of the modern development of structured ring spectra. The
multiplicative structure on the Koszul resolution gives us control over the
convergence of the spectral sequence and the multiplicative structure of
$\ell_*\ell$. In particular, it sheds light on the torsion.

The outline of the paper is as follows. We recall some basic facts about
complex cobordism, $MU$, in Section~\ref{sec:Recollections} and describe
the K\"unneth spectral sequence in Section~\ref{sec:KunnSS}. Some background
on the Bockstein spectral sequence is given in Section~\ref{sec:BocksteinSS}.
The multiplicative structure on the $\mathrm{E}^2$-term of this spectral
sequence is made precise in section \ref{sec:GenKoszul} where we introduce
the Koszul resolution we will use later in terms of its multiplicative
generators.  We study the torsion part in $\ell_*\ell$ and the torsion-free
part separately. The investigation of ordinary and $L$-homology of $\ell$ in
Section~\ref{sec:H*l} leads to the identification of the $p$-torsion in
$\ell_*\ell$ with the $u$-torsion where $\ell_* = \Z_{(p)}[u]$ with $u$ being
in degree $2p-2$. In Section~\ref{sec:ConnHomo} we show how to exploit the
cofibre sequence
$$
\ell \xrightarrow{p} \ell \lra \ell/p
$$
to analyse the K\"unneth spectral sequence and relate the simpler spectral
sequence for $\ell/p$ to that for $\ell$. To that end we prove an auxiliary
result on connecting homomorphisms in the K\"unneth spectral sequence, which
is analogous to the well-known geometric boundary theorem (see for
instance~\cite[chapter 2, \S 3]{Ravenel:greenbook}). We summarize our
calculation of $\ell_*\ell$ at the end of that section.

In Section~\ref{sec:ASS} we use classical tools from the Adams spectral
sequence in order to study torsion phenomena in $\ell_*\ell$. We use the fact
that the $p$ and $u$-torsion is all simple to show that the K\"unneth spectral
sequence for $\ell_*\ell$ collapses at the $\mathrm{E}^2$-term and that there
are no extension issues. We can describe the torsion in $\ell_*\ell$ in terms
of familiar elements which are certain coaction-primitives in the
$\F_p$-homology of $\ell$.

We summarize our results on the  multiplicative structure on
$\ell_*\ell$ at the end of Section~\ref{sec:l*l-Formulae}  
where we establish
congruence relations in the zero line of the K\"unneth spectral
sequence and  describe the map from the torsion-free part of $\ell_*\ell$
to $\ell_*\ell \otimes \mathbb{Q}$. Taking this together with the
explicit formul{\ae} of the multiplication in the torsion part in
$\ell_*\ell$ gives a rather comprehensive, though not complete,
description of the multiplicative structure of $\ell_*\ell$.

\section{Recollections on $MU$ and $\ell$}\label{sec:Recollections}

Throughout, we will assume all spectra are localized at~$p$ for some odd
prime~$p$.

Let $ku$ denote connective complex $K$-theory and let $\ell$ be the Adams
summand, also known as $BP\<1\>$, so that
\[
ku_{(p)}\sim\bigvee_{0\leq i\leq p-2}\Sigma^{2i}\ell.
\]
We have $\ell_*=\pi_*\ell=\Z_{(p)}[u]$ with $u\in\ell_{2(p-1)}$. We will
denote the Adams summand of $KU_{(p)}$ by $L$; then $L_*=\ell_*[u^{-1}]$.

Let us recall some standard facts for which convenient sources
are~\cite{JFA:Chicago,Wilson:Sampler}. Since $\ell$ is complex oriented,
\[
\ell_*MU=\ell_*[m'_n:n\geq1],
\]
where $m'_n\in\ell_{2n}MU$ agrees with the $m^\ell_n$ of Adams~\cite{JFA:Chicago}.
By the Hattori-Stong theorem, the Hurewicz homomorphism $MU_*\lra\ell_*MU$
is a split monomorphism, so we will view $MU_*$ as a subring of $\ell_*MU$.
Now
\[
MU_*=\Z_{(p)}[x_n:n\geq1],
\]
where $x_n\in MU_{2n}$ and using Milnor's criterion for polynomial generators
of $MU_*$ we can arrange that
\[
x_n\equiv
\begin{cases}
pm'_{p^k-1}\bmod{\mathrm{decomposables}}& \text{if $n=p^k-1$ for some $k$}, \\
\ph{p}m'_n\ph{-1}\bmod{\mathrm{decomposables}}& \text{otherwise}.
\end{cases}
\]
In fact, we can take $x_{p^k-1}=v_k$ to be the Hazewinkel generator which
lies in $BP_*\subset MU_*$. The following formula recursively determines
the Hurewicz image of~$v_k$ in $H_*MU=\Z_{(p)}[m_k:k\geq1]$:
\begin{equation}\label{eqn:Hazewinkel-MU}
v_k=pm_{p^k-1}-\sum_{1\leq j\leq k-1}m_{p^j-1}v_{k-j}^{p^j}.
\end{equation}
In $H_*BP$ with $\ell_k=m_{p^k-1}$, this corresponds to the familiar formula
\begin{equation}\label{eqn:Hazewinkel-BP}
v_k=p\ell_k-\sum_{1\leq j\leq k-1}\ell_jv_{k-j}^{p^j}.
\end{equation}

We note that
\begin{equation}\label{eqn:l*MUmodxn}
\ell_*MU/(x_n:\text{$n\neq p^k-1$ for any $k$})=\ell_*[t_k:k\geq1]=\ell_*BP,
\end{equation}
where $t_k\in\ell_{2p^k-2}BP$ is the image of the standard polynomial
generator $t_k\in BP_*BP$ of~\cite{JFA:Chicago}.

Now recall that the natural complex orientation of $\ell$ factors as
\[
\sigma\:MU\lra BP\lra\ell
\]
and we can choose the generators $x_n$ so that
\[
\sigma_*(x_n)=
\begin{cases}
u& \text{if $n=p-1$}, \\
0& \text{otherwise}.
\end{cases}
\]
In particular, the kernel of the map $BP_*\lra\ell_*$ is the ideal
generated by the Hazewinkel generators $v_2,v_3,\ldots$.

We can also find useful expressions for the $v_n$. Using standard
formul\ae{} for the right unit $\eta_R\:BP_*\lra BP_*BP$ which can
be found in~\cite{Wilson:Sampler}, we have for $n\geq2$,
\begin{equation}\label{eqn:l*BP-vn}
v_n=pt_n+ut_{n-1}^p-u^{p^{n-1}}t_{n-1}+ps'_n+us''_n,
\end{equation}
where $s'_n,s''_n\in\Z_{(p)}[u,t_1,\ldots,t_{n-1}]$. We also have
$v_1=pt_1+u$.

We now make some useful deductions. To ease notation we write $v_n$
for the image $\mathrm{e}(v_n)\in E_*BP$ of $v_n$ under the Hurewicz
homorphism $\mathrm{e}\:BP_*\lra E_*BP$.
\begin{prop}\label{prop:RegSeq-p}
In the ring $\Q\otimes\ell_*BP$, the sequence $v_2,v_3,\ldots,v_n,\ldots$
is regular and
\[
\Q\otimes\ell_*BP/(v_n:n\geq2)=\Q\otimes\ell_*[t_1]=\Q\otimes\ell_*[v_1].
\]
\end{prop}
\begin{proof}
For each $n\geq1$, $pt_n$ is a polynomial generator for
$\Q\otimes\ell_*BP=\Q\otimes\ell_*[t_i:i\geq1]$ over $\Q\otimes\ell_*$.
\end{proof}
\begin{prop}\label{prop:RegSeq-u}
In the ring $L_*BP$, the sequence $v_2,v_3,\ldots,v_n,\ldots$ is regular
and
\[
L_*BP/(v_n:n\geq2)=
L_*[t_k:k\geq1]/(t_{n}^p-u^{p^{n}-1}t_{n}
                  +pu^{-1}s'_{n+1}+s''_{n+1}+pu^{-1}t_{n+1}:n\geq1).
\]

In the ring $L_*BP/(p)$, the sequence $v_2,v_3,\ldots,v_n,\ldots$ is
regular and
\[
L_*BP/(p,v_n:n\geq2)
    =L_*/(p)[t_k:k\geq1]/(t_{n}^p-u^{p^{n}-1}t_{n}+s''_{n+1}:n\geq1).
\]
\end{prop}
\begin{proof}
The case of $L_*BP$ follows from Proposition~\ref{prop:Koszul-LandweberExact},
but here is a direct proof.

By induction, for each $m\geq2$ there is a monomorphism of $L_*$-algebras
\[
L_*BP/(v_n:m\geq n\geq2)\lra\Q\otimes L_*[t_1,t_k:k\geq m+1]
\]
under which the image of $v_{m+1}$ is clearly not a zero-divisor. This shows
that the $v_n$ form a regular sequence in $L_*BP$.
\end{proof}

\section{A K\"unneth spectral sequence for $\ell_*\ell$}
\label{sec:KunnSS}

We will describe a calculation of $\ell_*\ell=\pi_*(\ell\wedge\ell)$ that
makes use of the K\"unneth spectral sequence of~\cite{EKMM} for $MU$-modules.
This is different from the approach taken by Kane~\cite{Kane:AMSMem}, and we
feel it offers some insight into the form of answer, especially with regard
to multiplicative structure.

For any ring spectrum $E$ and $MU$-ring spectrum $F$, there is a homologically
graded spectral sequence
\begin{equation}\label{eqn:KunnSS}
\mathrm{E}^2_{s,t}=\Tor_{s,t}^{MU_*}(\pi_*(E\wedge MU),\pi_*F)
\Lra\pi_*((E\wedge MU)\wedge_{MU}F)\iso\pi_*(E\wedge F)=E_*F.
\end{equation}
Note that this spectral sequence is actually multiplicative~\cite{AB&AL:ASS}.
Taking $E=F=\ell$ we obtain a multiplicative spectral sequence
\begin{equation}\label{eqn:KunnSS-l*l}
\mathrm{E}^2_{s,t}=\Tor_{s,t}^{MU_*}(\pi_*(\ell\wedge MU),\pi_*\ell)
\Lra\ell_*\ell.
\end{equation}

Now consider the $MU_*$-module $\ell_*$. We can assume that the complex
orientation gives rise to a ring isomorphism
\[
MU_*/(x_n:n\neq p-1)\xrightarrow{\;\iso\;}\ell_*.
\]
There is a Koszul resolution of $\ell_*$ as a module over $MU_*$,
\[
\Lambda_{MU_*}(e_r:0<r\neq p-1)\lra\ell_*\ra0,
\]
where $\Lambda_{MU_*}(e_r:0<r\neq p-1)$ is the exterior algebra generated
by elements $e_r$ of bidegree $(1,2r)$ whose differential $d$ is the derivation
which satisfies $d(e_r)=x_r$.

For arbitrary $E$ and $F=\ell$, the $\mathrm{E}^2$-term of the spectral
sequence~\eqref{eqn:KunnSS} is the homology of the complex
\[
E_*MU\otimes_{MU_*}\Lambda_{MU_*}(e_r:0<r\neq p-1)\iso\Lambda_{E_*MU}(e_r:0<r\neq p-1)
\]
with differential $\id\otimes d$ which corresponds to the differential~$d$
taking values in the latter complex. From~\eqref{eqn:l*MUmodxn} we find
that the homology of this complex is
\begin{equation} \label{eqn:mutobp}
\mathrm{H}_*(\Lambda_{E_*MU}(e_r:0<r\neq p-1),d)=
\mathrm{H}_*(\Lambda_{E_*BP}(\epsilon_r:r\geq2),d),
\end{equation}
where $\epsilon_r$ has bidegree $(1,2p^r-2)$ and $d(\epsilon_r)=v_r$.
\begin{prop}\label{prop:Koszul-Exact}
Suppose that the $E$-theory Hurewicz images of the $v_k$ in $E_*BP$ with
$k\geq2$ form a regular sequence. Then the complex
\[
\Lambda_{E_*BP}(\epsilon_r:r\geq2)\lra E_*BP/(v_r:r\geq2)\ra0
\]
is acyclic and
\begin{equation}\label{eqn:Tor-Reg}
\Tor^{MU_*}_{s,*}(E_*MU,\ell_*)=
\begin{cases}
E_*BP/(v_r:r\geq2)& \text{\rm if $s=0$}, \\
\ph{abcabcabc}0 & \text{\rm otherwise}.
\end{cases}
\end{equation}
Therefore the K\"unneth spectral sequence of~\eqref{eqn:KunnSS}
degenerates to give an isomorphism
\[
E_*BP/(v_r:r\geq2)\xrightarrow{\;\iso\;}E_*\ell.
\]
\end{prop}

The regularity condition of this result occurs for the cases
$E=\ell\Q,\,L/p$ by Propositions~\ref{prop:RegSeq-u} and~\ref{prop:RegSeq-p}.
We do not have a proof that it holds for the case $E=L$, however
the following provides a substitute.
\begin{prop}\label{prop:Koszul-LandweberExact}
Suppose that $E$ is a $p$-local Landweber exact spectrum. Then the complex
\[
\Lambda_{E_*BP}(\epsilon_r:r\geq2) \lra E_*BP/(v_r:r\geq2)\ra0
\]
is acyclic and the conclusion of {\rm Proposition~\ref{prop:Koszul-Exact}}
is valid.
\end{prop}
\begin{proof}
The hypothesis means that the functor $E_*\otimes_{MU_*}(\ph{-})$ is exact
on the category of left $MU_*MU$-comodules. Then
\begin{align*}
E_*MU\otimes_{MU_*}\Lambda_{MU_*}(e_r:r\geq2)&\iso
E_*\otimes_{MU_*}MU_*MU\otimes_{MU_*}\Lambda_{MU_*}(e_r:r\geq2)  \\
&\iso E_*\otimes_{MU_*}\Lambda_{MU_*MU}(e_r:r\geq2) \\
&\iso E_*\otimes_{MU_*}\Lambda_{MU_*BP}(\epsilon_r:r\geq2).
\end{align*}
But now it is easy to see that the sequence $\eta_R(v_2),\eta_R(v_3),\ldots$
is regular in $MU_*BP$ since
\[
\eta_R(v_k)=v_k\bmod{(t_i:i\geq1)}
\]
and $v_2,v_3,\ldots$ is a regular sequence in $MU_*$. Therefore
\[
\Lambda_{MU_*BP}(\epsilon_r:r\geq2)\lra MU_*BP/(\eta_R(v_r):r\geq2)\ra0
\]
is exact and so is
\[
E_*\otimes_{MU_*}\Lambda_{MU_*BP}(\epsilon_r:r\geq2)
                  \lra E_*\otimes_{MU_*}MU_*BP/(\eta_R(v_r):r\geq2)\ra0,
\]
since in each homological degree, $\Lambda_{MU_*BP}(\epsilon_r:r\geq2)$
is a free $MU_*BP$-module and therefore an $MU_*MU$-comodule. From this
we obtain the result.
\end{proof}

Of course, this result applies when $E=L$. Later we will also consider
some cases where these regularity conditions do not hold.

\section{Bockstein spectral sequences}\label{sec:BocksteinSS}

We follow~\cite[p158]{Weibel} in this account. Let~$R$ be a graded commutative
ring and suppose that we have an exact couple of graded $R$-modules
\[
\xymatrix{
A^0_*\ar[rr]^{x\.}&&A^0_*\ar[dl]^{j^0} \\
&B^0_*\ar[ul]^{\delta^0}&
}
\]
where $\delta^0$ is a map of   degree $|x| - 1$ and $x\.$ is multiplication
by $x\in R$. Then there are inductively defined exact couples
\[
\xymatrix{
A^r_*\ar[rr]^{x\.}&&A^r_*\ar[dl]^{j^r} \\
&B^r_*\ar[ul]^{\delta^r}&
}
\]
and an associated spectral sequence $(B^r,d^r)$ with
$B^{r+1}_*=\mathrm{H}(B^r_*,d^r)$. For each $r\geq1$, there are exact
sequences
\begin{equation}\label{eqn:BSS-Er}
0\ra A^0_n/(xA^0_{n-|x|} + {}_{x^r}A^0_n)\xrightarrow{\bar{j^r}}B^r_n
            \xrightarrow{\delta^r} {}_{x}A^0_{n+|x|+r-1}
\cap x^rA^0_{n+|x|+r-r|x|-1} \ra0,
\end{equation}
where
\[
{}_{x^r}A^0_n=\ker(x^r\:A^0_{n} \lra A^0_{n+r|x|}),
\quad
{}_{x^\infty}A^0_n=\bigcup_{r\geq1}{}_{x^r}A^0_n.
\]
In particular, if $B^1_n=B^\infty_n=0$ for some $n$, we obtain the following:
\begin{align}
\label{eqn:BSS-Vanishing-torsion}
{}_{x^\infty}A_n&={}_{x}A_n,  \\
\label{eqn:BSS-Vanishing}
\ker\delta^0&=\ker d^0=\im j^0.
\end{align}

\section{Generalized Koszul complexes and Bockstein spectral sequences}
\label{sec:GenKoszul}

Let $R$ be a commutative ring and $x\in R$ a non-zero divisor which is
also not a unit. Let $w_1,w_2,w_3,\ldots$ be a (possibly finite) regular
sequence in~$R$ which reduces to a regular sequence in $R/(x)$.

The Koszul complex $(\Lambda_R(e_r:r\geq1),d)$ whose differential is the
$R$-derivation determined by $d(e_r)=w_r$ provides a resolution
\[
\Lambda_R(e_r:r\geq1)\lra R/(w_r:r\geq1)\ra0
\]
of $R/(w_r:r\geq1)$ by $R$-modules.

Now consider the sequence $xw_1,xw_2,xw_3,\ldots$ which is not regular
in~$R$ since for $s>r$,
\[
w_r(xw_s)=w_s(xw_r).
\]
The Koszul complex $(\Lambda_R(e'_r:r\geq1),d')$ with differential
satisfying $d'(e'_r)=xw_r$ is no longer exact but does augment onto
$R/(xw_r:r\geq1)$. Notice that there is a monomorphism of $R$-dga's
\[
j\:\Lambda_R(e'_r:r\geq1)\lra\Lambda_R(e_r:r\geq1);
\quad
j(e'_r)=xe_r,
\]
and this covers the reduction map $R/(xw_r:r\geq1)\lra R/(w_r:r\geq1)$.
Using this, we will view $\Lambda_R(e'_r:r\geq1)$ as a subcomplex of
$\Lambda_R(e_r:r\geq1)$. We want to determine the homology of
$(\Lambda_R(e'_r:r\geq1),d')$.

Suppose that $z\in\Lambda_R(e'_r:r\geq1)_n$ with $n>0$ and $d'(z)=0$.
Then working in $\Lambda_R(e_r:r\geq1)$ we have $d(j(z))=0$, so by
exactness of the latter complex, there is an element
\[
y=
\sum_{1\leq i_1<i_2<\cdots<i_{n+1}}
              y_{i_1,i_2,\ldots,i_{n+1}}e_{i_1}e_{i_2}\cdots e_{i_{n+1}}
\in\Lambda_R(e_r:r\geq1)_{n+1}
\]
for which $d(y)=j(z)$. But
\[
d(y)=
\sum_{\substack{1\leq i_1<i_2<\cdots<i_{n+1}\\1\leq k\leq n+1}}
(-1)^kw_{i_k}y_{i_1,i_2,\ldots,i_{n+1}}e_{i_1}e_{i_2}
                                 \cdots\hat{e}_{i_k} \cdots e_{i_{n+1}}.
\]
Since we have
\[
j(z)=\sum_{1\leq i_1<i_2<\cdots<i_n}
                  x^n z_{i_1,i_2,\ldots,i_n} e_{i_1}e_{i_2} \cdots e_{i_n},
\]
using the regularity assumption we find that each $y_{i_1,i_2,\ldots,i_{n+1}}$
has the form
\[
y_{i_1,i_2,\ldots,i_{n+1}}=x^n y'_{i_1,i_2,\ldots,i_{n+1}}
\]
for some $y'_{i_1,i_2,\ldots,i_{n+1}}\in R$ and therefore
\[
z=
\sum_{\substack{1\leq i_1<i_2<\cdots<i_{n+1}\\1\leq k\leq n+1}}
(-1)^kw_{i_k}y'_{i_1,i_2,\ldots,i_{n+1}}
             e'_{i_1}e'_{i_2}\cdots\hat{e'}_{i_k}\cdots e'_{i_{n+1}}.
\]
Notice that
\[
xz=
d'\(\sum_{1\leq i_1<i_2<\cdots<i_{n+1}}
             y'_{i_1,i_2,\ldots,i_{n+1}}e'_{i_1}e'_{i_2}\cdots e'_{i_{n+1}}\),
\]
hence $x$ annihilates the $n$-th homology of $\Lambda_R(e'_r:r\geq1)$,
hence it is an $R/(x)$-module spanned by the elements
\begin{equation}\label{eqn:Delta-Defn}
\Delta_x(i_1,i_2,\ldots,i_{n+1})=
\sum_{1\leq k\leq n+1}
           (-1)^kw_{i_k}e'_{i_1}e'_{i_2}\cdots\hat{e'}_{i_k}\cdots e'_{i_{n+1}}
\end{equation}
for collections of distinct integers $i_1,i_2,\ldots,i_{n+1}\geq1$.
Clearly for a permutation $\sigma\in\mathrm{S}_{n+1}$,
\[
\Delta_x(i_{\sigma(1)},i_{\sigma(2)},\ldots,i_{\sigma(n+1)})
                             = \sign\sigma\,\Delta_x(i_1,i_2,\ldots,i_{n+1}).
\]
Thus we will often restrict attention to indexing sequences satisfying
\[
1\leq i_1<i_2<\cdots<i_{n+1}.
\]
These elements satisfy some further additive and multiplicative relations.
\begin{prop}\label{prop:GenKosComp-Relations}
Let $r,s \geq 2$ and suppose that $i_1,i_2,\ldots,i_r\geq1$ and
$j_1,j_2,\ldots,j_s\geq1$ are sequences of distinct integers. Let
\[
t=\#\{i_1,i_2,\ldots,i_r\}\cup\{j_1,j_2,\ldots,j_s\}
\]
and write
\[
\{k_1,k_2,\ldots,k_t\}=\{i_1,i_2,\ldots,i_r\}\cup\{j_1,j_2,\ldots,j_s\}
\]
with $1\leq k_1<k_2<\cdots<k_t$. Then the following identities are satisfied
in each of $\Lambda_R(e'_r:r\geq1)$ and $\mathrm{H}_*(\Lambda_R(e'_r:r\geq1),d')$.
\begin{subequations}\label{eqn:GenKosComp-Relations}
\begin{multline}\label{eqn:GenKosComp-RelationsP}
\Delta_x(i_1,i_2,\ldots,i_r)\Delta_x(j_1,j_2,\ldots,j_s) = \\
\begin{cases}
0& \text{\rm if $t\leq r+s-2$}, \\[4pt]
(-1)^a w_{k_m}\Delta_x(k_1,k_2,\ldots,k_t)
                                    &\text{\rm if $t=r+s-1$ \& $k_m=i_a=j_b$},\\[4pt]
\ds\sum_{j=1}^r(-1)^{j+s+1}w_{i_j}
                \Delta_x(i_1,i_2,\ldots,\hat{i}_j,\ldots i_r,j_1,j_2,\ldots,j_s)
                                    & \text{\rm if $t=r+s$},
\end{cases}
\end{multline}
\begin{equation}\label{eqn:GenKosComp-RelationsS}
\sum_{j=1}^r(-1)^j w_{i_j}\Delta_x(i_1,i_2,\ldots,\hat{i}_j,\ldots i_r)=0.
\end{equation}
\end{subequations}
\end{prop}

\begin{thm}\label{thm:GenKosComp-Homology}
The homology of $(\Lambda_R(e'_r:r\geq1),d')$ is given by
\[
\mathrm{H}_n(\Lambda_R(e'_r:r\geq1),d')=
\begin{cases}
R/(xw_r:r\geq1)& \text{\rm if $n=0$}, \\[4pt]
R/(x)\{\Delta_x(i_1,i_2,\ldots,i_{n+1}):1\leq i_1<i_2<\cdots<i_{n+1}\}
                                            & \text{\rm if $n>0$}
\end{cases}
\]
where in the second case, the $R/(x)$-module is generated by the
$\Delta_x(i_1,i_2,\ldots,i_{n+1})$ indicated, subject to relations
given in~\eqref{eqn:GenKosComp-RelationsS}.
\end{thm}
\begin{proof}
Consider the long exact sequence obtained by taking homology of the
exact sequence
\[
0\ra R\otimes_R\Lambda_R(e'_r:r\geq1)\xrightarrow{}
              R\otimes_R\Lambda_R(e'_r:r\geq1)\xrightarrow{}
                    R/(x)\otimes_R\Lambda_R(e'_r:r\geq1)\ra0.
\]
The associated exact couple has
\begin{align*}
A^0_*&=\mathrm{H}_*(\Lambda_R(e'_r:r\geq1),d'), \\
B^0_*&=\mathrm{H}_*(\Lambda_{R/(x)}(e'_r:r\geq1),d')
=\Lambda_{R/(x)}(e'_r:r\geq1).
\end{align*}
Making use of the formula $d^0e'_r=w_r$ we find that
\[
B^1_*=R/(x,w_1,w_2,\ldots),
\]
and therefore the $x$-torsion in $A^0_*$ is all simple.
\end{proof}
Notice that the quotient $R$-module $R/(xw_r:r\geq1)$ has $x$-torsion,
as does the higher homology, at least if the sequence of $w_r$'s has
at least two terms.

\section{Ordinary and $L$-homology of $\ell$}\label{sec:H*l}

We can compute $H_*\ell$ using the spectral sequence
$(\mathrm{E}^r_{*,*}(H),d^r)$ obtained from~\eqref{eqn:KunnSS} by
taking $E=H=H\Z_{(p)}$ and $F=\ell$. This can be compared with the
spectral sequence $(\mathrm{E}^r_{*,*}(H\Q),d^r)$ for $H\Q_*\ell$
making use of the morphism of spectral sequences
\[
\mathrm{E}^r_{*,*}(H)\lra\mathrm{E}^r_{*,*}(H\Q)
\]
induced by the natural map $H\lra H\Q$. We will also consider the
spectral sequence $(\mathrm{E}^r_{*,*}(\bar{H}),d^r)$ associated
with $\bar{H}=H\F_p$.

By~\eqref{eqn:Hazewinkel-BP}, the sequence $v_2,v_3,\ldots,v_n,\ldots$
in the polynomial ring $H\Q_*BP=\Q[\ell_i:i\geq1]$ is regular. So by
Proposition~\ref{prop:Koszul-Exact} we have
\begin{equation} \label{eq:rational}
\mathrm{E}^2_{s,*}(H\Q)=
\begin{cases}
\Q[\ell_i:i\geq1]/(v_k:k\geq2) &\text{if $s=0$}, \\
0&\text{otherwise}.
\end{cases}
\end{equation}
Hence this spectral sequence collapses at $\mathrm{E}^2$ and we have
\[
H\Q_*\ell=\Q[\ell_1]=\Q[v_1],
\]
where $v_1=p\ell_1$. The image of $\ell_n$ in $H\Q_*\ell$ can be
recursively computed with the aid of the following formula derived
from~\eqref{eqn:Hazewinkel-BP}:
\begin{equation}\label{eqn:tn-Recursion}
\ell_n=\frac{v_1^{p^{n-1}}\ell_{n-1}}{p}.
\end{equation}
So we have
\begin{equation}\label{eqn:H*l-ln}
\ell_n=\frac{v_1^{(p^n-1)/(p-1)}}{p^n}
=p^{p^{n-1}+p^{n-2}+\cdots+p+1-n}\,\ell_1^{(p^n-1)/(p-1)}.
\end{equation}
Notice that for a monomial in the $\ell_j$'s in $H\Q_{2m(p-1)}\ell$,
we have
\[
\ell_1^{r_1}\cdots\ell_n^{r_n}=\frac{v_1^m}{p^{r_1+2r_2+\cdots+nr_n}},
\]
for which
\[
r_1+2r_2+\cdots+nr_n\leq
r_1+r_2\frac{p^2-1}{p-1}+\cdots+r_n\frac{p^n-1}{p-1}=m.
\]
This calculation shows that the images of the monomials in the $\ell_j$'s
in $H\Q_{2m(p-1)}\ell$ are contained in the cyclic $\Z_{(p)}$-module
generated by $\ell_1^m=v_1^m/p^m$. Turning to the spectral sequence
$\mathrm{E}^r_{*,*}(H)$, we see that
\[
\mathrm{E}^2_{0,*}(H)=H_*BP/(v_j:j\geq2)
\]
and the natural map
\[
H_{2m(p-1)}BP/(v_j:j\geq2)\lra H\Q_{2m(p-1)}BP/(v_j:j\geq2)
\]
has image equal $\Z_{(p)}\ell_1^m$. In~\cite{JFA:Chicago}, the analogous
result for $ku$ was obtained using the Adams spectral sequence.
\begin{prop}\label{prop:H*l->HQl}
For $m\geq0$,
\[
\im[H_{2m(p-1)}\ell\lra H\Q_{2m(p-1)}\ell]=\Z_{(p)}\ell_1^m=\Z_{(p)}\frac{v_1^m}{p^m}.
\]
Hence,
\[
\im[H_*\ell\lra H\Q_*\ell]=\Z_{(p)}[\ell_1]=\Z_{(p)}[v_1/p].
\]
\end{prop}

The spectral sequence $(\mathrm{E}^r_{*,*}(\bar{H}),d^r)$ is easy to
determine. As $v_k=0$ in $\bar{H}_*BP$, we find that
\[
\mathrm{E}^\infty_{*,*}(\bar{H})=\mathrm{E}^2_{*,*}(\bar{H})
                            =\Lambda_{\bar{H}_*BP}(\epsilon_r:r\geq2).
\]
Thus we recover the well-known result that
\[
\bar{H}_*\ell=\F_p[t_k:k\geq1]\otimes_{\F_p}\Lambda_{\F_p}(\epsilon_r:r\geq2),
\]
where $t_k$ has degree $2p^k-2$ and $\epsilon_r$ has degree $2p^r-1$.


From Propositions~\ref{prop:RegSeq-u} and~\ref{prop:Koszul-Exact}
we have
\begin{align*}
\Tor^{MU_*}_{*,*}(L_*MU,\ell_*)&=L_*BP/(v_r:r\geq2), \\
\Tor^{MU_*}_{*,*}(\bar{L}_*MU,\ell_*)&=\bar{L}_*BP/(v_r:r\geq2),
\end{align*}
where $\bar{L}=L/p$ denotes the spectrum $L$ smashed with the mod~$p$
Moore spectrum. As a consequence, the K\"unneth spectral sequences
for $L_*\ell$ and $\bar{L}_*\ell$ degenerate to give
\[
L_*BP/(v_r:r\geq2)\iso L_*\ell,
\quad
\bar{L}_*BP/(v_r:r\geq2)\iso\bar{L}_*\ell.
\]
Since $L_*MU$ is a free $\Z_{(p)}$-module, multiplication by~$p$
gives an exact sequence of right $MU_*$-modules
\[
0\ra L_*MU\xrightarrow{p}L_*MU\lra\bar{L}_*MU\ra0
\]
which in turn induces a long exact sequence on the homological
functor $\Tor^{MU_*}_*(\ph{X},\ell_*)$ which collapses to the
short exact sequence
\[
0\ra\Tor^{MU_*}_{0,*}(L_*MU,\ell_*)\xrightarrow{\;p\;}
            \Tor^{MU_*}_{0,*}(L_*MU,\ell_*)\lra
                 \Tor^{MU_*}_{0,*}(\bar{L}_*MU,\ell_*)\ra0.
\]
From this we see that there is a short exact sequence
\[
0\ra L_*\ell\xrightarrow{\;p\;}L_*\ell\lra\bar{L}_*\ell\ra0.
\]
On tensoring with $\Q$ we easily see that
$\Q\otimes\ell_*\ell\lra\Q\otimes L_*\ell$\/ is a monomorphism.
Hence we have
\begin{prop}\label{prop:L*l-NoTorsion}
$L_*\ell$ has no $p$-torsion and the natural map\/
$\ell_*\ell\lra L_*\ell$\/ induces an exact sequence
\[
0\ra{}_{p^\infty}(\ell_*\ell)\lra\ell_*\ell\lra L_*\ell.
\]
\end{prop}
\begin{cor}\label{cor:L*l-NoTorsion}
We have
\[
{}_{p^\infty}(\ell_*\ell)={}_{u^\infty}(\ell_*\ell).
\]
\end{cor}
\begin{proof}
Since $\ell_*\lra L_*=\ell_*[u^{-1}]$ is a localization,
we have $L_*\ell=\ell_*\ell[u^{-1}]$ and
\[
\ker(\ell_*\ell\lra L_*\ell)={}_{u^\infty}(\ell_*\ell),
\]
hence ${}_{u^\infty}(\ell_*\ell)={}_{p^\infty}(\ell_*\ell)$.
\end{proof}

\section{Connecting homomorphisms in the K\"unneth spectral sequence}
\label{sec:ConnHomo}

In order to gain control over the $p$-torsion in
$\Tor_{*,*}^{MU_*}(\ell_*MU, \ell_*)$, we will exploit the cofibre
sequence
\begin{equation}\label{eqn:ell/p-seq}
\ell \xrightarrow{p} \ell \xrightarrow{\varrho} \bar{\ell}
\xrightarrow{\delta} \Sigma\ell,
\end{equation}
where $\bar{\ell} = \ell/p$. To this end we will relate the geometric
connecting morphisms of cofibre sequences to  morphisms of K\"unneth
spectral sequences. The method of proof we use in this part is analogous
to that of the geometric boundary theorem in~\cite[II.3]{Ravenel:greenbook}.

Let $W$ be a cofibrant $R$-module which we fix from now on. Then for
any $R$-module $Z$ there is a K\"unneth spectral sequence with
\[
\mathrm{E}^2_{s,t}(Z)=\Tor_{s,t}^{R_*}(Z_*,W_*)
                    \Longrightarrow \pi_*(Z\wedge_R W).
\]
\begin{lem} \label{lem:connfree}
Let
$$
X \xrightarrow{f} Y \xrightarrow{g} Z \xrightarrow{h} \Sigma X
$$
be a cofibre sequence of $R$-modules with $X \simeq \bigvee_{i=1}^m \Sigma^{n_i} R$
and $\pi_*f$ surjective. Then there is a map of K\"unneth spectral sequences
$$
\mathrm{E}^r_{s,t}(Y) \xrightarrow{\psi^r}\mathrm{E}^r_{s-1,t}(\Sigma^{-1}Z)
\quad (r \geq 2),
$$
such that $\psi^2$ is the connecting homomorphism
\[
\Tor_{s,t}^{R_*}(Y_*, W_*) \xrightarrow{\;\cong\;}
                          \Tor_{s-1,t}^{R_*}({(\Sigma^{-1}Z)}_*, W_*).
\]
\end{lem}
\begin{proof}
Since $\pi_*f$ is surjective, there is a short exact sequence
$$
0 \ra {(\Sigma^{-1}Z)}_* \lra \bigoplus_{i=1}^m \Sigma^{n_i}R_* \lra Y_* \ra 0.
$$
This induces a long exact sequence of $\Tor$-groups, in which every third
term is trivial, because $\bigoplus_{i=1}^m \Sigma^{n_i}R_*$ is $R_*$-free.
Therefore we have an isomorphism
\[
\Tor_{s,t}^{R_*}(Y_*, W_*) \xrightarrow{\;\cong\;}
                           \Tor_{s-1,t}^{R_*}({(\Sigma^{-1}Z)}_*, W_*).
\]
On the level of projective resolutions, we can splice a resolution
$P_{\bullet,*}$ for $Y_*$ together with a resolution $Q_{\bullet,*}$ of
${(\Sigma^{-1}Z)}_*$ to obtain a trivial split resolution for
$\bigoplus_{i=1}^m \Sigma^{n_i}R_*$. Thus we obtain a map between exact
couples and so obtain the desired map of spectral sequences.
\end{proof}
\begin{thm}\label{thm:connkuenneth}
Let
$$
X \xrightarrow{f} Y \xrightarrow{g} Z \xrightarrow{h} \Sigma X
$$
be a cofibre sequence with $\pi_*f$ surjective. Then there is an induced
map of K\"unneth spectral sequences
$$
\mathrm{E}^r_{s,t}(Y) \xrightarrow{\phi^r}\mathrm{E}^r_{s-1,t}(\Sigma^{-1}Z)
\quad (r \geq 2)
$$
such that $\phi^2$ is the connecting homomorphism
\[
\Tor_{s,t}^{R_*}(Y_*, W_*) \xrightarrow{\;\cong\;}
                            \Tor_{s-1,t}^{R_*}({(\Sigma^{-1}Z)}_*, W_*).
\]
\end{thm}
\begin{proof}
Choose a map $f'\: \bigvee_{i=1}^m \Sigma^{n_i} R \lra Y$ with $\pi_*f'$
surjective and consider the cofibre sequence
$$
\bigvee_{i=1}^m \Sigma^{n_i} R \xrightarrow{f'} Y \xrightarrow{j} \cone(f').
$$
By Lemma~\ref{lem:connfree} there is a map of K\"unneth spectral sequences
$$
\mathrm{E}^r_{s,t}(Y) \xrightarrow{\psi^r}
\mathrm{E}^r_{s-1,t}(\Sigma^{-1}\cone(f')).
$$
As $\pi_*f$ is surjective, the composition $g \circ f'$ is trivial and there
is a factorization $g = \xi \circ j$.
\[
\xymatrix{
{} & {\cone(f')} \ar@{.>}[dr]^{\xi}& {} & {} \\
{X} \ar[r]^{f} & {Y} \ar[r]^{g} \ar[u]^{j} & {Z}\ar[r]^{h} & {\Sigma X} \\
{} & {\bigvee_{i=1}^m \Sigma^{n_i} R} \ar[u]^{f'} & {} & {}
}
\]
Now we may define $\phi^r$ to be $(\Sigma^{-1}\xi)_* \circ \psi^r$.
\end{proof}

For the connective Adams summand $\ell$, we will consider the cofibre
sequence
\[
\ell \xrightarrow{\varrho} \bar{\ell} \xrightarrow{\delta} \Sigma\ell
                                     \xrightarrow{\Sigma p} \Sigma\ell
\]
obtained from~\eqref{eqn:ell/p-seq}. The reduction map $\varrho$ is surjective
in homotopy and therefore we can apply Theorem~\ref{thm:connkuenneth} to
obtain a map of K\"unneth spectral sequences
\[
\mathrm{E}^r_{s,t}(\bar{\ell}\wedge MU) \xrightarrow{\phi^r}
                                  \mathrm{E}^r_{s-1,t}(\ell\wedge MU)
\quad(r \geq 2).
\]
In particular, this yields a connecting homomorphism
\[
\phi^2\: \Tor_{s,t}^{MU_*}(\bar{\ell}_*MU, \ell_*) \lra
                                      \Tor_{s-1,t}^{MU_*}(\ell_*MU, \ell_*).
\]
\begin{thm}\label{thm:ptorsion}
Each $p$-torsion element of\/ $\Tor_{*,*}^{MU_*}(\ell_*MU, \ell_*)$ is
the image of an element of \/ $\Tor_{*+1,*}^{MU_*}(\bar{\ell}_*MU,
\ell_*)$  under the
connecting homomorphism and is an infinite cycle.
\end{thm}
\begin{proof}
Making use of the long exact sequence on $\Tor$-groups associated to the
short exact sequence
\[
0 \ra \ell_* \xrightarrow{p} \ell_* \xrightarrow{\varrho_*} \bar{\ell}_* \ra 0,
\]
the claim about the $p$-torsion in
$\Tor_{*,*}^{MU_*}(\ell_*MU, \ell_*)$
follows.

We will prove that the elements $\Delta_u(i_1,\ldots,i_n)$ are infinite
cycles in the K\"unneth spectral sequence for $\bar{\ell}_*\ell$. Our
proof is by induction on $n\geq2$. For $n=2$ the elements $\Delta_u(i_1,i_2)$
are infinite cycles for degree reasons. Now suppose that for all $n\leq k$,
the  $\Delta_u(i_1,\ldots,i_n)$ are infinite cycles. From
Proposition~\ref{prop:GenKosComp-Relations} we know that
$$
w_{i_2}\Delta_u(i_1,\ldots,i_{k+1}) =
\Delta_u(i_1,i_2)\Delta_u(i_2,\ldots,i_{k+1}).
$$
By assumption, the two factors are infinite cycles and therefore their
product is an infinite cycle as well. The scalar factor $w_{i_2}$ acts
as a regular element on the $R/(u)$-module generated by the $\Delta_u$
elements, so we can conclude that $\Delta_u(i_1,\ldots,i_{k+1})$ has
to be an infinite cycle as well.

So the K\"unneth spectral sequence for $\bar{\ell}_*\ell$ collapses at
the $\mathrm{E}^2$-page and as the connecting  homomorphism induces a
map of spectral sequences, every $p$-torsion class in
$\Tor_{*,*}^{MU_*}(\ell_*MU, \ell_*)$ has to be an infinite cycle.
\end{proof}

Theorem~\ref{thm:ptorsion} gives an explicit description of the
$p$-torsion classes in $\Tor_{*,*}^{MU_*}(\ell_*MU, \ell_*)$ as the
image of the elements $\Delta_u(i_1,\ldots,i_n)$ in
$\Tor_{*,*}^{MU_*}(\bar{\ell}_*MU, \ell_*)$ under the boundary homorphism.

\begin{cor}\label{cor:ptorsion}
Since the $\Delta_u(i_1,\ldots,i_n)$ generate the K\"unneth spectral
sequence for $\bar{\ell}_*\ell$ additively, this spectral sequence
collapses at the $\mathrm{E}^2$-page. For the same reasons, the
K\"unneth spectral sequence for $\ell_*\ell$ collapses as well.
\end{cor}
\begin{rem} \label{rem:collapsing}
To summarize, the calculation of the rational homology of $\ell$ in
\eqref{eq:rational} tells us that the torsion-free part of $\ell_*\ell$
has to have its origin in the zero-line of the K\"unneth spectral sequence.
The torsion part is imported from the K\"unneth spectral sequence for
$\bar{\ell}_*\ell$ via the geometric boundary result. The K\"unneth spectral
sequence for $\ell_*\ell$ collapses at the $\mathrm{E}^2$-page. Furthermore,
Corollary~\ref{cor:simpletorsion} implies that there are no extension problems.
\end{rem}

%

\section{Detecting homotopy in the Adams spectral sequence}\label{sec:ASS}

In this section we recall some results about the classical Adams spectral
sequence for $\ell_*\ell$. We make heavy use of standard facts about Hopf
algebras and the Steenrod algebra~\cite{Milnor&Moore,Milnor}. In the following
we generically write $I$ for identity morphisms, $\phi$ for products and
actions, $\psi$ for coproducts and coactions, $\eta$ for units and $\epsilon$
for counits and we use $\bar{x}$ for the antipode on an element $x$.
Undecorated tensor products are taken over the ground field.

We write $\bar{H}_*(\ )$ for $H_*(\ ;\F_p)$ and $\StA_*$ for the dual Steenrod
algebra,
\[
\StA_*=\F_p[\zeta_n:n\geq1]\otimes\Lambda(\bar\tau_n:n\geq0),
\]
where the coaction is given by
\[
\psi(\zeta_n)=\sum_{i=0}^n\zeta_i\otimes \zeta_{n-i}^{p^i},
\quad \psi(\bar\tau_n)=
   1\otimes\bar\tau_n+\sum_{i=0}^n\bar\tau_i \otimes \zeta_{n-i}^{p^i}.
\]
The sub-comodule algebra
\[
\StB_*=\F_p[\zeta_n:n\geq1]\otimes\Lambda(\bar\tau_n:n\geq2)
\]
gives rise to a quotient Hopf algebra
\[
\StE_*=\StA_*/\!/\StB_*=\Lambda(\alpha,\beta),
\]
where $\alpha,\beta$ are the residue classes of $\bar\tau_0,\bar\tau_1$
respectively. Then
\[
\StB_*=\StA_*\square_{\StE_*}\F_p.
\]
Now the natural map $\ell\lra \bar{H}$ induces an isomorphism
\[
\bar{H}_*(\ell) \xrightarrow{\;\iso\;} \StB_* \subset \StA_*
\]
and there are isomorphisms of $\StA_*$-comodule algebras
\begin{equation}\label{eqn:ComIsos}
\bar{H}_*(\ell\wedge\ell) \xrightarrow{\;\iso\;} \bar{H}_*(\ell) \otimes \bar{H}_*(\ell)
\xrightarrow{\;\iso\;} \StB_* \otimes \StB_*
\xrightarrow{\;\iso\;} \StA_* \square_{\StE_*} \StB_*.
\end{equation}
The $\mathrm{E}_2$-term of the Adams spectral sequence converging
to $\pi_*(\ell\wedge\ell)=\ell_*\ell$ has the form
\[
\mathrm{E}_{s,t}^2=\Cotor^{\StA_*}_{s,t}(\F_p,\bar{H}_*(\ell\wedge\ell))
\iso \Cotor^{\StA_*}_{s,t}(\F_p,\StA_* \square_{\StE_*} \StB_*)
\]
and so by making use of a standard change of rings result, we have
\begin{equation}\label{eqn:ll-ASSE2}
\mathrm{E}_{s,t}^2 \iso \Cotor^{\StE_*}_{s,t}(\F_p, \StB_*).
\end{equation}

Note that by results of~\cite{Kane:AMSMem}, the torsion in $\ell_*\ell$
is detected by the edge homomorphism (which is essentially the Hurewicz
homomorphism) into the $0$-line
\[
\mathrm{E}_{0,*}^2 \iso \Cotor^{\StE_*}_{0,*}(\F_p, \StB_*)
                                   = \F_p \square_{\StE_*} \StB_*.
\]
The map involved here is obtained by composing the following $\StA_*$-comodule
algebra homomorphisms and suitably restricting the codomain:
\begin{multline*}
\pi_*(\ell\wedge\ell) \lra \bar{H}_*(\ell) \otimes \bar{H}_*(\ell)
\xrightarrow{\;\iso\;} \StB_* \otimes \StB_*
\xrightarrow{I \otimes \psi} \StB_* \otimes (\StA_* \square_{\StE_*} \StB_*)  \\
\xrightarrow{\phi \otimes I}
\StA_* \square_{\StE_*} \StB_*
\xrightarrow{\ph{I \otimes \phi}}
\StE_* \square_{\StE_*} \StB_*
\xrightarrow{\;\iso\;} \StB_*.
\end{multline*}
The final isomorphism is the composition
\[
\StE_* \square_{\StE_*} \StB_* \xrightarrow{\;\incl\;} \StE_* \otimes \StB_*
\xrightarrow{\epsilon\otimes I} \F_p \otimes \StB_*
\xrightarrow{\;\iso\;} \StB_*.
\]
A careful check of what the composition does on primitives shows that it
can be expressed as
\begin{equation} \label{eq:themap}
\pi_*(\ell \wedge \ell) \lra \bar{H}_*(\ell \wedge \ell)
\xrightarrow{(\nu \wedge \mathrm{id})_*} \bar{H}_*(\ell),
\end{equation}
where $\nu\: \bar{H} \wedge \ell \lra \bar{H}$ is the natural pairing. In particular,
this implies that the image of the Hurewicz map for
$\ell \wedge \ell$ maps monomorphically into $\bar{H}_*(\ell)$.

It will be useful to know how to compute the inverse of the map
\[
\F_p\square_{\StA_*} (\StB_* \otimes \StB_*) \lra \F_p\square_{\StE_*}\StB_*.
\]
This is just
\[
\F_p\square_{\StE_*}\StB_* \xrightarrow{\incl} \F_p \otimes \StB_*
         \xrightarrow{I \otimes \psi} \F_p \otimes (\StA_* \otimes \StB_*),
\]
whose image is in fact contained in $\F_p \square_{\StA_*} (\StB_* \otimes \StB_*)$.

Given these results, we can use them to detect elements of $\ell_*\ell$
in $\StB_*$, in particular we can detect the torsion this way. To do this,
we need to understand $\StB_*$ as an $\StE_*$-comodule, in particular the
non-trivial $\StE_*$-parallelograms of the form
\begin{equation}\label{eqn:StE*-parallelogram}
\xymatrix{
&&&&&x\ar[dllll]_{-\beta}\ar[dl]^{\alpha}  \\
&x''\ar[dl]_{\alpha}&&&x'\ar[dllll]^{\beta}&  \\
x'''&&&&&
}
\end{equation}
in which the $\StE_*$-coaction satisfies
\begin{equation} 
\psi(x)=
1\otimes x + \alpha\otimes x' - \beta\otimes x''
                            + \beta\alpha\otimes x''',
\quad
\psi(x')=1\otimes x' + \beta\otimes x''',
\quad
\psi(x'')=1\otimes x'' + \alpha\otimes x'''.
\end{equation}
Then $x'''$ is an element of $\F_p\square_{\StE_*}\StB_*$ which corresponds
to an $H\F_p$ wedge summand in $\ell \wedge \ell$ and a correponding torsion
element. Of course, these elements can be expressed in terms of the homology
action of $Q_0$ and $Q_1$, \ie,
\[
x'=Q_0x,\quad x''=-Q_1x,\quad x'''=Q_1Q_0x.
\]

Now by Margolis~\cite[chapter~18 theorem~5]{Margolis-book} dualized to a
homology version for $\StE_*$-comodules tells us that $\StB_*$ uniquely
decomposes into a coproduct of comodules isomorphic (up to grading) to
$\StE_*$, together with a comodule containing no free summand and isomorphic
to a coproduct of lightning flash comodules. The latter summand does not
concern us for now since all the torsion in $\ell_*\ell$ comes from the
$H\F_p$ wedge summands as above corresponding to the free summand. In fact,
Adams and Priddy~\cite[proof of proposition~3.12]{JFA&SP:BSU} determine
the stable type of the lightning flash comodules, in particular, the stable
class of the $\StE_*$-comodule $\StB_*$ is shown to be
\begin{equation}\label{eqn:StB-stabletype}
\bigotimes_{r\geq0}(1 + L_r +L_r^2 + \cdots + L_r^{p-1}),
\end{equation}
where
\[
L_r = \Sigma^{a(r)}J^{b(r)},
\quad
a(r)+b(r)=2(p-1)p^r,
\quad
b(r) = p^{r-1} + \cdots +p+1.
\]
Here $J = \StE_*/\F_p$ is the coaugmentation coideal of $\StE_*$,
represented by the following diagram
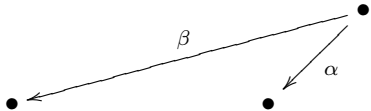
\[
\xymatrix{
&&&&\bullet\ar[dllll]_{\beta}\ar[dl]^{\alpha} \\
\bullet&&&\bullet&
}
\]
and $\Sigma$ is the trivial comodule $\F_p$ assigned degree~$1$. Furthermore,
all products are tensor products over $\F_p$ taken in the stable comodule
category of 
$\StE_*$.

Now the most obvious candidates for the tops of $\StE_*$-parallelograms
are the elements
\[
\bar\tau_{i_1} \bar\tau_{i_2} \cdots \bar\tau_{i_{n+1}}
\qquad
(1 < i_1< i_2<\cdots < i_{n+1},\quad n\geq1).
\]
These can be multiplied by monomials in the $\zeta_j$ to obtain others.
\begin{thm}\label{thm:H*l-Parallelograms}
Consider the $\F_p$-vector subspace $\mathcal{V}\subset
\F_p\square_{\StE_*}\StB_*$
spanned by  $\F_p[\zeta_i:i\geq1]$-scalar multiples of the
elements $1$ and
\begin{equation}\label{eqn:H*l-Parallelograms}
Q_1Q_0(\bar\tau_{i_1} \bar\tau_{i_2}  \cdots \bar\tau_{i_{n+1}})
\qquad
(1 < i_1< i_2<\cdots < i_{n+1},\quad n\geq1).
\end{equation}
Then $\mathcal{V}$ consists of all the elements in $\F_p\square_{\StE_*}\StB_*$
which are the images of torsion elements under the composition of the
Hurewicz homomorphism
$\pi_*(\ell \wedge \ell)\lra \bar{H}_*(\ell \wedge \ell)$ and the
identification of the homology $\bar{H}_*(\ell \wedge \ell)$ with
$\F_p\square_{\StE_*}\StB_*$.
\end{thm}
\begin{proof}
Clearly $\F_p[\zeta_i:i\geq1] \subset \F_p\square_{\StE_*}\StB_*$.
Now we know that the K\"unneth spectral sequence for $\ell_*\ell$ collapses
and there are no additive extension problems. We need to understand the
mod~$p$ Hurewicz images of elements represented by the elements arising
from the $\Delta_u(i_1,\ldots,i_{s+2})$ in $\Tor^{MU_*}_{s+1,*}(\bar{\ell}_*MU,\ell_*)$,
since these will give an additive basis for the $p$-torsion in $\ell_*\ell$.
$$
\xymatrix{
{\Tor_{s+1,*}^{MU_*}(\bar{\ell}_*MU, \ell_*)} \ar@/^/[dr]_{\delta} & & & \\
& {\Tor_{s,*}^{MU_*}({\ell}_*MU, \ell_*)} \ar[d] \ar@{.}[r]&
{\pi_*(\ell \wedge \ell)} \ar@/^/[rd]
\ar[d] & \\
& {\Tor_{s,*}^{MU_*}((\bar{H} \wedge \ell)_*MU, \ell_*)} \ar[d] \ar@{.}[r]&
{\pi_*(\bar{H} \wedge \ell \wedge \ell)} \ar[d] & {\F_p \square_{\StA_*}
\bar{H}_*(\ell \wedge \ell)} \ar[l] \ar[d]\\
& {\Tor_{s,*}^{MU_*}(\bar{H}_*MU, \ell_*)} \ar@{.}[r] & {\bar{H}_*(\ell)}
& {\F_p \square_{\StE_*} \bar{H}_*(\ell)} \ar[l] }
$$
The K\"unneth spectral sequence~\eqref{eqn:KunnSS} for $E_*\ell$ is natural
for maps of ring spectra $E \lra F$. Therefore the map~\eqref{eq:themap}
corresponds in the spectral sequence to the composition of the two vertical
maps in the left column in the diagram above. As the Hurewicz homomorphism
has its image in the primitives of $\bar{H}_*(\ell \wedge \ell)$, it follows
that the elements $\Delta_u(i_1, \ldots, i_{s+2})$ up to a unit  map to
\[
Q_0Q_1(\bar{\tau}_{i_1} \cdots \bar{\tau}_{i_{s+2}})
= \sum_{1 < t < r \leq s+2}
(-1)^{r+t} (\zeta_{i_r}\zeta_{i_t-1}^p-\zeta_{i_t}\zeta_{i_r-1}^p) \,
\bar\tau_{i_1} \bar\tau_{i_2} \cdots \hat{\bar\tau}_{i_t}
              \cdots \hat{\bar\tau}_{i_r} \cdots \bar\tau_{i_{s+2}}.
\qedhere
\]
\end{proof}
\begin{rem}
Since the torsion in $\pi_*(\ell \wedge \ell)$ maps injectively into
$\F_p \square_{\StA_*} (\StB_* \otimes \StB_*)$ which in turn is
identified with $\F_p \square_{\StE_*} \StB_*$, Theorem
\ref{thm:H*l-Parallelograms} shows that the elements
$Q_1Q_0(\bar\tau_{i_1} \bar\tau_{i_2}  \cdots \bar\tau_{i_{n}})$ with
$n \geq 3$ correspond to nilpotent elements; only elements of the form
$Q_1Q_0(\bar\tau_{r}\bar\tau_{s})$ are not nilpotent.
\end{rem}
From Corollary~\ref{cor:L*l-NoTorsion} we know that the $p$-torsion and
$u$-torsion in $\ell_*\ell$ agree. We recall a fact
from~\cite[proposition 9.1]{Kane:AMSMem}.
\begin{prop} \label{prop:simpletorsion}
All torsion in $\ell_*\ell$ is simple, \ie, for every torsion-class
$x \in \ell_*\ell$ we have $px = 0$ which is equivalent to $ux = 0$.
\end{prop}
\begin{cor} \label{cor:simpletorsion}
The K\"unneth spectral sequence for $\ell_*\ell$ collapses at the
$\mathrm{E}^2$-page and there are no non-trivial extensions.
\end{cor}
\begin{ex}\label{ex:first-torsion}
For every prime $p$, the first torsion class in
$\Tor_{*,*}^{BP_*}(\ell_*BP, \ell_*)$ occurs in degree $2(p^3 + p^2 - p -1)$
and this class survives to $\ell_*\ell$.
The lowest degree element appearing as the bottom of a parallelogram is
\[
Q_1Q_0(\bar\tau_2\bar\tau_3)=\zeta_2^{p+1} -\zeta_1^p\zeta_3.
\]
The coaction map $\psi$ sends this element to the Hurewicz image of
the corresponding torsion element of $\ell_*\ell$ in
$\bar{H}_*(\ell \wedge \ell)$.
\end{ex}

\section{Multiplicative structure of $\ell_*\ell$}\label{sec:l*l-Formulae}

In this section we establish
congruence relations in the zero line of the K\"unneth spectral
sequence. These are derived in $BP_*BP$ and mapped under the natural
map. In  fact
they are first produced in $\Q\otimes BP_*BP$ then interpreted in the
subring $BP_*BP$.

We  describe the map from the torsion-free part of $\ell_*\ell$
to $\ell_*\ell \otimes \mathbb{Q}$ and summarize our results about the
multiplicative structure of $\ell_*\ell$ at the end of this section. 

It will be useful to have the following generalization of a well-known
result (which corresponds to the case where $t=1$).
\begin{lem}\label{lem:pku->p(k+1)u}
Let $R$ be a commutative ring, $p$ a prime and $t\in R$. If $x,y,z\in R$
satisfy $z\equiv px+ty\bmod(pt)$, then for all $k\geq0$,
\[
z^{p^k}\equiv p^{p^k}x^{p^k}+t^{p^k}y^{p^k}\bmod(p^{k+1}t).
\]
\end{lem}
\begin{proof}
We prove this by induction on $k$, the case $k=0$ being known. Suppose
it is true for some $k\geq0$. Choose a $w\in R$ for which
\[
z^{p^k}=p^{p^k}x^{p^k}+t^{p^k}y^{p^k}+p^{k+1}tw.
\]
Then working $\bmod(p^{k+2}t)$ we have
\begin{align*}
z^{p^{k+1}}&=(p^{p^k}x^{p^k}+t^{p^k}y^{p^k})^p+p^{k+2}t^pw^p
 + \sum_{1\leq i\leq p-1}
           \binom{p}{i}(p^{p^k}x^{p^k}+t^{p^k}y^{p^k})^{p-i}p^{k+1+i}t^iw^i  \\
& \equiv(p^{p^k}x^{p^k}+t^{p^k}y^{p^k})^p \\
& \equiv p^{p^{k+1}}x^{p^{k+1}}+t^{p^{k+1}}y^{p^{k+1}}
+ \sum_{1\leq i\leq p-1}\binom{p}{i}p^{ip^k}x^{ip^k}t^{(p-i)p^k}y^{(p-i)p^k} \\
& \equiv p^{p^{k+1}}x^{p^{k+1}}+t^{p^{k+1}}y^{p^{k+1}}.
\end{align*}
Hence the result holds for $k+1$.
\end{proof}

We will work with the Hazewinkel generators~$v_n$ of~\eqref{eqn:Hazewinkel-BP}.
The following standard formula for the right unit
$\eta_R\:\Q\otimes BP_*\lra\Q\otimes BP_*BP$ can be found in~\cite[p24]{Wilson:Sampler}:
\begin{equation}\label{eqn:etaR(ln)}
\eta_R(\ell_n)=\sum_{0\leq j\leq n}\ell_jt_{n-j}^{p^j}.
\end{equation}
On combining this with~\eqref{eqn:Hazewinkel-BP} we obtain
\[
\eta_R(v_n)=
\sum_{0\leq i\leq n}p\ell_it_{n-i}^{p^i}-
\sum_{\substack{1\leq i\leq n-1\\0\leq j\leq i}}\ell_jt_{i-j}^{p^j}\eta_R(v_{n-i})^{p^i}
\]
and hence
\begin{equation}\label{eqn:Vn-Recursion}
\eta_R(v_n) = \sum_{0\leq i\leq n}p\ell_it_{n-i}^{p^i} -
\sum_{0\leq i\leq n-1}\ell_it_{n-1-i}^{p^i}\eta_R(v_1)^{p^{n-1}} -
\sum_{\substack{1\leq i\leq n-2\\0\leq j\leq i}}\ell_jt_{i-j}^{p^j}\eta_R(v_{n-i})^{p^i}.
\end{equation}
\begin{rem}\label{rem:Vn-Recursion}
The left hand side of equation~\eqref{eqn:Vn-Recursion} lies in
$BP_*BP\subset\Q\otimes BP_*BP$, therefore so does the right hand side.
However, because of the presence of denominators in the terms involving
the $\ell_r$, care needs to be exercised when using this equation. For
example, since $\cp_r=p^r\ell_r\in BP_*$ we can certainly deduce that in
$BP_*BP$ modulo the ideal $(\eta_R(v_2),\ldots,\eta_R(v_{n-1}))\ideal BP_*BP$,
\begin{multline*}
p^{n-1}\eta_R(v_n)\equiv \\
\sum_{0\leq i\leq n}p^{n-i}\cp_it_{n-i}^{p^i}
-\sum_{0\leq i\leq n-1}p^{n-1-i}\cp_it_{n-1-i}^{p^i}\eta_R(v_{1})^{p^{n-1}}
\bmod{(\eta_R(v_2),\ldots,\eta_R(v_{n-1}))}.
\end{multline*}
We will see later that similar phenomena in $\ell_*BP$ give rise
to congruences in $\ell_*\ell$.
\end{rem}

We will now derive some formul\ae{} in $\ell_*BP$. The natural map
of ring spectra $BP\lra\ell$ is determined on homotopy by
\begin{equation}\label{eqn:l*-vr}
v_r\longmapsto
\begin{cases}
u& \text{if $r=1$}, \\
0& \text{otherwise}.
\end{cases}
\end{equation}
Recalling~\eqref{eqn:H*l-ln}, we see that in $\im[H_*\ell\lra H\Q_*\ell]$,
the logarithm series for the factor of $\ell$ is
\[
\log^\ell T=\sum_{n\geq0}\ell_nT^{p^n}
=\sum_{n\geq0}\frac{u^{p^{n-1} + \cdots + p + 1}}{p^n}T^{p^n}.
\]
We can project~\eqref{eqn:Vn-Recursion} into $\ell_*BP$, with $\eta_R$
being replaced by the $\ell$-theory Hurewicz homomorphism
$\underline{\ell}\:BP_*\lra\ell_*BP$. This yields
\begin{multline*}
\underline{\ell}(v_n)=
pt_n-t_{n-1}\underline{\ell}(v_{1})^{p^{n-1}} \\
+\sum_{1\leq i\leq n}\frac{u^{p^{i-1}+\cdots+p+1}t_{n-i}^{p^i}}{p^{i-1}}
-\sum_{1\leq i\leq n-1}
  \frac{u^{p^{i-1}+\cdots+p+1}t_{n-1-i}^{p^i}\underline{\ell}(v_1)^{p^{n-1}}}{p^{i}} \\
-\sum_{1\leq i\leq n-2}t_{i}\underline{\ell}(v_{n-i})^{p^i}
-\sum_{\substack{1\leq i\leq n-2\\1\leq j\leq i}}
\frac{u^{p^{j-1}+\cdots+p+1}t_{i-j}^{p^j}\underline{\ell}(v_{n-i})^{p^i}}{p^j}.
\end{multline*}
and the equivalent formula
\begin{multline}\label{eqn:l*BP-Vn-Recursion-2}
\underline{\ell}(v_n)=
pt_n+(ut_{n-1}^p-\underline{\ell}(v_{1})^{p^{n-1}}t_{n-1}) \\
+\sum_{1\leq i\leq n-1}
  \frac{u^{p^{i-1}+\cdots+p+1}(u^{p^{i}}t_{n-1-i}^{p^{i+1}}
             -\underline{\ell}(v_1)^{p^{n-1}}t_{n-1-i}^{p^i})}{p^{i}} \\
-\sum_{1\leq i\leq n-2}t_{i}\underline{\ell}(v_{n-i})^{p^i}
-\sum_{\substack{1\leq i\leq n-2\\1\leq j\leq i}}
\frac{u^{p^{j-1}+\cdots+p+1}t_{i-j}^{p^j}\underline{\ell}(v_{n-i})^{p^i}}{p^j}.
\end{multline}
Thus we have
\begin{align*}
\underline{\ell}(v_2)&=
pt_2+(ut_{1}^p-\underline{\ell}(v_{1})^{p}t_{1})
                      +\frac{u(u^{p}-\underline{\ell}(v_1)^{p})}{p} \\
&=pt_2+(1-p^{p-1})ut_{1}^p-\underline{\ell}(v_{1})^{p}t_{1}
                       -\sum_{1\leq i\leq p-1}\binom{p}{i}p^{i-1}u^{p+1-i}t_1^i.
\end{align*}

By the Hattori-Stong theorem, the element $\underline{\ell}(v_n)\in\ell_*BP$
is not divisible by~$p$, but notice that on multiplying by $p^{n-2}$ we have
\begin{multline*}
p^{n-2}\underline{\ell}(v_n)=
p^{n-1}t_n+p^{n-2}(ut_{n-1}^p-\underline{\ell}(v_{1})^{p^{n-1}}t_{n-1}) \\
+\sum_{1\leq i\leq n-1}
  \frac{u^{p^{i-1}+\cdots+p+1}(u^{p^{i}}t_{n-1-i}^{p^{i+1}}
             -\underline{\ell}(v_1)^{p^{n-1}}t_{n-1-i}^{p^i})}{p^{i-n+2}} \\
-\sum_{1\leq i\leq n-2}p^{n-2}t_{i}\underline{\ell}(v_{n-i})^{p^i}
-\sum_{\substack{1\leq i\leq n-2\\1\leq j\leq i}}
        p^{n-2-j}u^{p^{j-1}+\cdots+p+1}t_{i-j}^{p^j}\underline{\ell}(v_{n-i})^{p^i}.
\end{multline*}
and so
\begin{multline*}
p^{n-1}t_n+p^{n-2}(ut_{n-1}^p-\underline{\ell}(v_{1})^{p^{n-1}}t_{n-1}) \\
+\sum_{1\leq i\leq n-1}
 \frac{u^{p^{i-1}+\cdots+p+1}(u^{p^{i}}t_{n-1-i}^{p^{i+1}}
             -\underline{\ell}(v_1)^{p^{n-1}}t_{n-1-i}^{p^i})}{p^{i-n+2}}
           \equiv0
          \quad\bmod{(\underline{\ell}(v_2),\ldots,\underline{\ell}(v_n))}.
\end{multline*}

Using the identity $\underline{\ell}(v_1)=u+pt_1$ and the resulting
congruences (see Lemma~\ref{lem:pku->p(k+1)u}),
\[
\underline{\ell}(v_1)^{p^m}\equiv u^{p^m}\bmod{(p^{m+1})}\qquad(m\geq1),
\]
we deduce that when $n\geq2$,
\begin{multline}\label{eqn:l*BP-Vn-Recursion-congruence}
\underline{\ell}(v_n)\equiv
(pt_n-p^{p^{n-1}}t_1^{p^{n-1}}t_{n-1})+(ut_{n-1}^p-u^{p^{n-1}}t_{n-1}) \\
+\sum_{1\leq i\leq n-2}
\frac{u^{p^{i-1}+\cdots+p+1}(u^{p^i}t_{n-1-i}^{p^{i+1}}-u^{p^{n-1}}t_{n-1-i}^{p^i})}{p^{i}}
\ph{\sum_{1\leq i\leq n-2}v_{n-i}}   \\
-\sum_{1\leq i\leq n-2}t_{i}\underline{\ell}(v_{n-i})^{p^i}
-\sum_{\substack{1\leq i\leq n-2\\1\leq j\leq i}}
\frac{u^{p^{j-1}+\cdots+p+1}t_{i-j}^{p^j}\underline{\ell}(v_{n-i})^{p^i}}{p^j}
\quad\bmod{(pu)}.
\end{multline}
Thus when $n=2$ we have
\begin{align*}
\underline{\ell}(v_2)&\equiv(pt_2-p^{p}t_1^{p}t_{1})+(ut_{1}^p-u^{p}t_{1})
\quad\bmod{(pu)} \\
&\equiv ut_{1}^p-u^{p}t_{1}\quad\bmod{(p)}.
\end{align*}

When working in the image of the rationalization map
$H_*(\ell\wedge\ell)\lra H\Q_*(\ell\wedge\ell)$, we will denote by~$u$ and~$v$
the images of $u\in\ell_{2p-2}$ under the left and right units for $\ell\wedge\ell$.

Now reinterpreting~\eqref{eqn:Vn-Recursion} in $H\Q_*(\ell\wedge\ell)$,
for each $n\geq2$ we have $\eta_R(v_n)\longmapsto 0$ and so
\begin{multline*}
pt_n+ut_{n-1}^p+
\sum_{1\leq h\leq n-1}\frac{u^{(p^h+p^{h-1}+\cdots+p+1)}t_{n-h-1}^{p^{h+1}}}{p^h} \\
=t_{n-1}v^{p^{n-1}}+
\sum_{1\leq k\leq n-1}
        \frac{u^{(p^{k-1}+p^{k-2}+\cdots+p+1)}t_{n-1-k}^{p^{k}}v^{p^{n-1}}}{p^{k}}.
\end{multline*}
On rearranging this, we obtain
\begin{equation}\label{eqn:l*l-Units}
pt_n=v^{p^{n-1}}t_{n-1}-ut_{n-1}^p
+\sum_{1\leq k\leq n-1}\frac{u^{(p^{k-1}+p^{k-2}+\cdots+p+1)}
(v^{p^{n-1}}t_{n-1-k}^{p^{k}}-u^{p^k}t_{n-k-1}^{p^{k+1}})}{p^{k}}.
\end{equation}
For small values of $n=1$ we have
\begin{align*}
pt_1&=v-u, \\
pt_2&=v^pt_1-ut_1^p+\frac{u(v^p-u^p)}{p}, \\
pt_3&=v^{p^2}t_2-ut_2^p+\frac{u(v^{p^2}t_1^{p}-u^pt_1^{p^2})}{p}
                             +\frac{u^{p+1}(v^{p^2}-u^{p^2})}{p^2},  \\
pt_4&=v^{p^3}t_3-ut_3^p+\frac{u(v^{p^3}t_2^{p}-u^pt_2^{p^2})}{p}
  +\frac{u^{p+1}(v^{p^3}t_1^{p^2}-u^{p^2}t_1^{p^3})}{p^2}
          +\frac{u^{p^2+p+1}(v^{p^3}-u^{p^3})}{p^3}.
\end{align*}

We want to draw some general conclusions about these expressions.
\begin{lem}\label{lem:pu}
In $\ell_*\ell$, for $n\geq1$, we have the congruences
\begin{align}
pt_n &\equiv v^{p^{n-1}}t_{n-1}-ut_{n-1}^p\bmod{(pu)},
\label{eqn:pu-1}\\
pt_n-p^{p^{n-1}}t_1^{p^{n-1}}&\equiv u^{p^{n-1}}t_{n-1}-ut_{n-1}^p\bmod{(pu)}.
\label{eqn:pu-2}
\end{align}
\end{lem}
\begin{proof}
We will prove this by induction on~$n$, the case $n=1$ being noted
above. So suppose that
\[
pt_k\equiv v^{p^{k-1}}t_{k-1}-ut_{k-1}^p\bmod{(pu)}.
\]
whenever~$1\leq k<n$ for some~$n$. Then for every such~$k$ we have
\[
v^{p^{k-1}}t_{k-1}\equiv ut_{k-1}^p\bmod{(p)}.
\]
By Lemma~\ref{lem:pku->p(k+1)u}, for every $m\geq1$,
\begin{align*}
(v^{p^{k-1}}t_{k-1})^{p^m}&\equiv(ut_{k-1}^p)^{p^m}\bmod{(p^{m+1})}, \\
\intertext{\ie,}
v^{p^{m+k-1}}t_{k-1}^{p^m}&\equiv u^{p^m}t_{k-1}^{p^{m+1}}\bmod{(p^{m+1})}.
\end{align*}
Now when $1\leq k\leq n-1$,
\[
v^{p^{n-1}}t_{n-1-k}^{p^{k}}-u^{p^k}t_{n-k-1}^{p^{k+1}}\equiv0\bmod{(p^{k+1})},
\]
hence in the formula for $pt_n$ in~\eqref{eqn:l*l-Units}, the summand
\[
u^{(p^{k-1}+p^{k-2}+\cdots+p+1)}
\frac{(v^{p^{n-1}}t_{n-1-k}^{p^{k}}-u^{p^k}t_{n-k-1}^{p^{k+1}})}{p^{k}}
\]
must be divisible by $pu$. Therefore we have the congruence
\[
pt_n\equiv v^{p^{n-1}}t_{n-1}-ut_{n-1}^p\bmod{(pu)}.
\]
Using the expansion
\[
v^{p^{n-1}}=u^{p^{n-1}}+\sum_{1\leq j\leq p^{n-1}}\binom{p^{n-1}}{j}u^{p^{n-1}-j}p^jt_1^j
\]
we obtain
\[
pt_n-p^{p^{n-1}}t_1^{p^{n-1}}\equiv u^{p^{n-1}}t_{n-1}-ut_{n-1}^p\bmod{(pu)}.
\qedhere
\]
\end{proof}
\subsection*{Summary}
Kane \cite[(19:6:1)]{Kane:AMSMem}, using Adams' criterion  
\cite[III,17.6]{JFA:Chicago},  worked out what the image of the
torsion-free part of
$\ell_*\ell$ is when we pass to $\ell_*\ell \otimes \mathbb{Q}$.
The generators for the image of $\ell_*\ell$ in  
$\ell_*\ell \otimes \mathbb{Q}$ are
$$ t_{n,i} = \frac{u^i v(v -(p-1)u)\cdot \ldots \cdot 
(v - (n-1)(p-1)u)}{p^i}, \quad 0 \leq i \leq \nu_p(n!).$$
Obviously, the relation $ut_{n,i} = pt_{n,i+1}$ holds and it is clear how to multiply elements like that.

To summarize our results on the multiplicative structure of $\ell_*\ell$, we have the following:
\begin{itemize}
    \item
    When we start with two non-torsion elements in $\ell_*\ell$, we
    can take their image in $\ell_*\ell \otimes \mathbb{Q}$,
    take their product there and interpret the result as a
    non-torsion element in $\ell_*\ell$.
    \item
    Any two elements coming from the zero-line of the K\"unneth
    spectral sequence multiply according to the congruence
    relations we specified in \eqref{eqn:l*-vr} up to \eqref{lem:pu}. These
    element might be torsion or non-torsion, but there is no
    non-torsion in higher filtrations.
    \item
    Torsion elements in non-zero filtration have their origin in
    the generators $\Delta_u$ and for these we spelled out the
    multiplication in \eqref{eqn:GenKosComp-RelationsP}.
    \item
    As the $\Delta_u$-expression allow coefficients from
    $\ell_*BP$, the multiplication of non-torsion elements in the
    zero-line with torsion elements in higher filtration is
    determined as well.
\end{itemize}
We agree that the recursive nature of the congruences for
$\ell_*BP \otimes_{BP_*} \ell_*$ might hamper the calculation, but
our approach leads to more information about the multiplication in
$\ell_*\ell$ than the known sources (compare \eg,~\cite[p.~76]{Kane:AMSMem}).

\end{document}